\documentclass[colorlinks]{mpi2015-cscpreprint}
\usepackage[american]{babel}
\usepackage{graphicx}
\usepackage{caption}
\usepackage{booktabs}
\usepackage{url}
\usepackage{xurl}
\usepackage[section]{placeins}
\usepackage{amssymb,amsthm,mathtools}
\usepackage{textcomp}
\usepackage{algorithm,algpseudocode}
\AtBeginDocument{%
  \theoremstyle{plain}
  \newtheorem{theorem}{Theorem}[section]
  
  \newtheorem{proposition}[theorem]{Proposition}
  \newtheorem*{conjecture}{Conjecture}
  \theoremstyle{definition}
  
  \theoremstyle{remark}
  
}
\makeatletter
\providecommand*{\theHALG@line}{}
\renewcommand*{\theHALG@line}{\thealgorithm.\arabic{ALG@line}}
\AtBeginDocument{%
  \g@addto@macro\ALG@beginalgorithmic{\def\@currentcounter{ALG@line}}%
  \crefalias{ALG@line}{line}%
  \crefrangeformat{line}{lines~#3#1#4--#5#2#6}}
\makeatother
\algnewcommand{\IIf}[1]{\State\algorithmicif\ #1\ \algorithmicthen}
\algnewcommand{\EndIIf}{\unskip\ \algorithmicend\ \algorithmicif}
\newcommand{\R}{\mathbb R}
\newcommand{\diag}{\operatorname{diag}}
\newcommand{\rank}{\operatorname{rank}}
\newcommand{\range}{\operatorname{range}}
\newcommand{\vol}{\operatorname{vol}}
\newcommand{\F}{\mathrm F}

\begin{document}
\title{Selecting \texorpdfstring{$2 \times 2$}{2x2} submatrices with optimally bounded inverses}
\author[$\ast$]{Bonan Sun}
\affil[$\ast$]{Max Planck Institute for Dynamics of Complex Technical Systems, Sandtorstr.~1, D-39106 Magdeburg, Germany.\authorcr%
  \email{bsun@mpi-magdeburg.mpg.de}, \orcid{0009-0002-2761-0661}}
\shorttitle{Selecting $2 \times 2$ submatrices with optimally bounded inverses}
\shortauthor{B. Sun}
\shortdate{}
\abstract{%
  For a matrix $U \in \R^{n \times k}$ ($n \ge k$) with orthonormal columns, the Goreinov--Tyrtyshnikov--Zamarashkin (GTZ) conjecture asserts that there exists a $k \times k$ submatrix $U_I$ with $\|U_I^{-1}\|_2 \le \sqrt{n}$.
  The conjecture remains open for $3 \le k \le n - 3$.
  The existence of such a submatrix in the case $k = 2$ was proved for $n=4$ by Nesterenko (2023) and recently for general $n$ by Sengupta and Pautov (2026).
  In this paper, we give a short, elementary, constructive proof of a weighted strengthening of this result, which yields a single-pass streaming algorithm that selects a row pair satisfying the $\sqrt{n}$ bound in $\mathcal O(n)$ operations and $\mathcal O(1)$ auxiliary memory.
  We further give an explicit matrix representation of all equality cases for $k = 2$ characterized by Nesterenko (2026) and use it to construct a single matrix family on which standard selection strategies (including global and local maximum-volume selection, Bischof--Stewart pivoted QR, column-pivoted QR, and LU with partial pivoting) simultaneously fail to attain the $\sqrt{n}$ bound.}
\novelty{%
   For matrices with two orthonormal columns, we prove a weighted row-selection bound by elementary row elimination, turn the proof into a single-pass $\mathcal O(n)$ streaming algorithm attaining the sharp $\sqrt{n}$ bound, and use a three-block representation of the equality cases to construct a single matrix family on which maximum-volume and standard greedy pivoting methods all exceed $\sqrt{n}$ by a constant factor.}
\maketitle

\section{Introduction}
\label{sec:intro}

Let $U \in \R^{n \times k}$ with $n \ge k \ge 1$ have orthonormal columns, $U^\top U = I_k$, and write $U_I \in \R^{k \times k}$ for the row submatrix indexed by $I \subseteq \{1,\ldots,n\}$ with $|I| = k$.
Selecting a ``good'' submatrix $U_I$ with small inverse spectral norm $\|U_I^{-1}\|_2$ is a classical problem in numerical linear algebra, with close connections to discrete empirical interpolation (DEIM~\cite{ChaturantabutSorensen2010} and QDEIM~\cite{DrmacGugercin2016}), strong rank-revealing factorizations~\cite{HongPan1992,GuEisenstat1996,Damle2026}, column subset selection~\cite{ChanHansen1992,GuEisenstat1996,BoutsidisMahoneyDrineas2009,BoutsidisDrineasMagdonIsmail2014,Osinsky2025,CortinovisKressner2026}, and (pseudo)skeleton, cross, or CUR approximations~\cite{GTZ1997theory,OsinskyZamarashkin2018,MikhalevOseledets2018}.
See~\cite{Simons2026} for a recent overview of open questions on index selection and rank-revealing factorizations.

Classical selection strategies provide quite different worst-case bounds on $\|U_I^{-1}\|_2$.
Standard sequential pivoting schemes such as column-pivoted QR (QRCP~\cite{BusingerGolub1965,DrmacGugercin2016}) and LU with partial pivoting (LUPP~\cite{ChaturantabutSorensen2010,DrmacGugercin2016}) can suffer from exponential growth in $k$ on Kahan-type matrices~\cite{GuEisenstat1996,ChaturantabutSorensen2010,DrmacGugercin2016}.
In contrast, Bischof--Stewart pivoted QR (BSQR), which dates back to 1990~\cite{Bischof1990,Stewart1990} and selects rows sequentially by greedily minimizing $\|U_{I_t}^\dagger\|_\F^2$ one row at a time ($U_{I_t}^\dagger$ here denotes the Moore--Penrose pseudoinverse), was recently shown to attain the bound $\sqrt{k(n-k)+1}$~\cite{Osinsky2025,Damle2026}.
The same bound can be achieved by choosing $I$ globally to maximize the volume $|\det U_I|$ (exact maxvol)~\cite{Knuth1985,GuEisenstat1996,GTZ1997theory,GTZ1997maxvol}.
Although exact maxvol is NP-hard~\cite{CivrilMagdonIsmail2009} for general $k$, greedily increasing $|\det U_I|$ via single-row swaps (as in strong rank-revealing QR~\cite{GuEisenstat1996} and the maxvol algorithm~\cite{Goreinov2010GoodSubmatrix}) attains $\sqrt{1 + (1+\varepsilon)^2 k(n-k)}$ for any swap tolerance $\varepsilon > 0$.
The bound $\sqrt{k(n-k)+1}$ is sharp for maximum-volume submatrices~\cite{GTZ1997theory}, and even the slightly smaller bounds obtained via interlacing polynomials~\cite{Xu2025} still scale as $\sqrt{k(n-k)}$.
Goreinov, Tyrtyshnikov, and Zamarashkin conjectured in~\cite{GTZ1997theory} that every $U$ in fact admits a submatrix with inverse spectral norm at most $\sqrt{n}$, independent of $k$.

\begin{conjecture}[Goreinov--Tyrtyshnikov--Zamarashkin (GTZ) conjecture]
For every $n \ge k \ge 1$ and every $U \in \R^{n \times k}$ satisfying $U^\top U = I_k$, there exists an index set $I \subseteq \{1,\ldots,n\}$ with $|I| = k$ such that $U_I$ is invertible and
\begin{equation}\label{eq:gtz}
 \|U_I^{-1}\|_2 \le \sqrt{n}.
\end{equation}
\end{conjecture}

The case $k = 1$ is trivial, as is $n = k$.
For $k = 2$, the classical bound $\sqrt{2n-3}$ strictly exceeds $\sqrt{n}$ for all $n \ge 4$.
The GTZ conjecture for $k = 2$ was proved for $n = 4$ by Nesterenko~\cite{Nesterenko2024} and recently settled for all $n \ge 2$ by Sengupta and Pautov~\cite{SenguptaPautov2026}, with sharpness for every $n \ge 3$ and the equality criterion established in~\cite{Nesterenko2026equality}.
By orthogonal-complement duality ($k \leftrightarrow n - k$)~\cite{GTZ1997theory}, \eqref{eq:gtz} also holds for $n = k + 1$ and $n = k + 2$.
Consequently, the conjecture remains open for $3 \le k \le n - 3$.
Interestingly, the conjecture fails in the complex case, where no bound of the form $C\sqrt{n}$ with a universal constant $C$ holds (see \href{https://github.com/ajt60gaibb/OpenProblemsInNLA/tree/main/randomized-and-low-rank-approximation/RA-18}{Problem~RA-18} in~\cite{Townsend2026openproblems}).

However, the proofs of the GTZ conjecture for $k=2$ in~\cite{Nesterenko2024,SenguptaPautov2026} establish only the existence of a row pair satisfying \eqref{eq:gtz} and do not provide an efficient selection algorithm.
In particular, the general-$n$ proof in~\cite{SenguptaPautov2026} relies on a contradiction argument via the Perron--Frobenius theorem, leaving open whether a $2 \times 2$ submatrix satisfying \eqref{eq:gtz} can be found in $\mathcal O(n)$ operations rather than by $\mathcal O(n^2)$ exhaustive search over all $\binom{n}{2}$ pairs.

In this paper, we answer this question through a constructive proof of the following weighted strengthening of \eqref{eq:gtz} for $k=2$.

\begin{theorem}\label{thm:weighted}
Let $n \ge 2$, let $U \in \R^{n \times 2}$ satisfy $U^\top U = I_2$, and let $w_1,\ldots,w_n \in (0, 1)$ satisfy $\sum_{r=1}^n w_r \le 1$.
Then there exists a two-element subset $I = \{i, j\} \subseteq \{1,\ldots,n\}$ ($i < j$) such that $U_I$ is invertible and
\begin{equation}\label{eq:weighted-bound}
 \bigl\|U_I^{-1}\diag(\sqrt{w_i},\sqrt{w_j})\bigr\|_2 \le 1.
\end{equation}
\end{theorem}

Setting $w_1 = \cdots = w_n = 1/n$ in \Cref{thm:weighted} recovers the GTZ bound \eqref{eq:gtz} for $k = 2$.
In \Cref{sec:proof}, we give a short, elementary proof of \Cref{thm:weighted} by direct row elimination (a Carath\'eodory-type reduction~\cite[Section~17]{Rockafellar1970}).
Allowing general positive weights keeps the problem closed under weight updates, which is crucial for turning the proof into a single-pass streaming algorithm in \Cref{sec:algorithm}.
In \Cref{sec:counterexamples}, we give an explicit matrix representation of all equality cases of \eqref{eq:gtz} for $k=2$ characterized in~\cite{Nesterenko2026equality}, in terms of three-block orthonormal bases obtained by repeating and rescaling the weighted three-row equality cases of \Cref{thm:weighted}.
Using three-block bases, we also construct a single matrix family on which exact maxvol, maxvol, QRCP, BSQR, and LUPP all exceed $\sqrt{n}$ as $n \to \infty$, whereas our streaming algorithm attains the global minimum and reduces $\|U_I^{-1}\|_2$ by an asymptotic factor of $(1 + \sqrt{5})/2$ relative to all five methods.
Finally, \Cref{sec:experiments} reports numerical experiments.

\section{Proof of \texorpdfstring{\Cref{thm:weighted}}{Theorem 1.1}}
\label{sec:proof}

\begin{proof}
Define $A = \diag(w_1^{-1/2},\ldots,w_n^{-1/2}) U \in \R^{n \times 2}$.
For any $I = \{i, j\}$ with $1 \le i < j \le n$, we have $U_I = \diag(\sqrt{w_i},\sqrt{w_j}) A_I$, so \eqref{eq:weighted-bound} holds if and only if $A_I$ is invertible with $\|A_I^{-1}\|_2 \le 1$.
We proceed by induction on $n \ge 2$.

The claim for $n = 2$ is trivial.
For $n = 3$, recall that a matrix $M \in \R^{m \times m}$ is invertible with $\|M^{-1}\|_2 \le 1$ if and only if $\|M x\|_2 \ge \|x\|_2$ for all $x \in \R^m$.
For any $x \in \R^2$, setting $z = (z_1, z_2, z_3)^\top = U x \in \R^3$ gives $\|x\|_2^2 = \|U x\|_2^2 = z_1^2 + z_2^2 + z_3^2$ and $\|A_{\{1,2\}} x\|_2^2 = z_1^2 / w_1 + z_2^2 / w_2$, so $A_{\{1,2\}}$ is invertible with $\|A_{\{1,2\}}^{-1}\|_2 \le 1$ if and only if
\[
 z_3^2 \le \frac{1 - w_1}{w_1} z_1^2 + \frac{1 - w_2}{w_2} z_2^2 \qquad \text{for all } x \in \R^2.
\]
Let $v = (v_1, v_2, v_3)^\top \in \range(U)^\perp$ be a unit vector.
Since $v^\top z = 0$ gives $v_3 z_3 = -(v_1 z_1 + v_2 z_2)$, the Cauchy--Schwarz inequality yields
\[
 v_3^2 z_3^2 = (v_1 z_1 + v_2 z_2)^2 \le \left(\frac{w_1 v_1^2}{1 - w_1} + \frac{w_2 v_2^2}{1 - w_2}\right)\left(\frac{1 - w_1}{w_1} z_1^2 + \frac{1 - w_2}{w_2} z_2^2\right).
\]
Hence $A_{\{1,2\}}$ is invertible with $\|A_{\{1,2\}}^{-1}\|_2 \le 1$ whenever $v_3 \ne 0$ and $w_1 v_1^2 / (1 - w_1) + w_2 v_2^2 / (1 - w_2) \le v_3^2$.
By permuting the rows of $U$ and $w$, we may assume without loss of generality that $v_1^2 / (1 - w_1) \le v_2^2 / (1 - w_2) \le v_3^2 / (1 - w_3)$.
Since $v$ is a unit vector, we have $v_3 \ne 0$, and $w_1 + w_2 \le 1 - w_3$ gives
\begin{equation}\label{eq:n3-key}
 \frac{w_1 v_1^2}{1 - w_1} + \frac{w_2 v_2^2}{1 - w_2} \le (w_1 + w_2)\,\frac{v_3^2}{1 - w_3} \le v_3^2,
\end{equation}
which completes the proof for $n = 3$.

For $n \ge 4$, assume that the claim holds for all row counts smaller than $n$, and let $a_r^\top$ denote the $r$-th row of $A$.
Since $U = \diag(\sqrt{w_1},\ldots,\sqrt{w_n}) A$, the condition $U^\top U = I_2$ reads $\sum_{r=1}^n w_r a_r a_r^\top = I_2$.
To eliminate at least one row while preserving this identity and the rows of $A$, observe that the $n \ge 4$ matrices $a_1 a_1^\top,\ldots,a_n a_n^\top$ lie in the three-dimensional space of $2 \times 2$ real symmetric matrices and are therefore linearly dependent.
Thus there exists a nonzero vector $\gamma = (\gamma_1,\ldots,\gamma_n)^\top \in \R^n$ satisfying $\sum_{r=1}^n \gamma_r a_r a_r^\top = 0$.
Without loss of generality, assume $\sum_{r=1}^n \gamma_r \le 0$, so at least one entry satisfies $\gamma_r < 0$.
Increasing $\theta \ge 0$ in $w_r + \theta \gamma_r$ until the first weight reaches zero, set
\[
 \theta_* = \min_{\gamma_r < 0} \frac{w_r}{-\gamma_r} > 0, \qquad w_r' = w_r + \theta_* \gamma_r \ge 0, \qquad r = 1,\ldots,n,
\]
and let $S = \{r \in \{1,\ldots,n\} : w_r' > 0\}$.
Since at least one $w_r'$ vanishes, by permuting the rows of $U$ and $w$ we may assume without loss of generality that $S = \{1,\ldots,m\}$ with $m \le n - 1$, and set $U' = \diag(\sqrt{w_1'},\ldots,\sqrt{w_m'}) A_S \in \R^{m \times 2}$.
To apply the induction hypothesis to $U'$ with weights $w_1',\ldots,w_m'$, we verify that $U'$ has orthonormal columns (so that $m \ge 2$), $w_r' \in (0, 1)$ for $r = 1,\ldots,m$, and $\sum_{r=1}^m w_r' \le 1$.
Indeed, we have
\[
 (U')^\top U' = \sum_{r=1}^m w_r' a_r a_r^\top = \sum_{r=1}^n (w_r + \theta_* \gamma_r) a_r a_r^\top = I_2,
\]
and $\sum_{r=1}^m w_r' = \sum_{r=1}^n w_r + \theta_* \sum_{r=1}^n \gamma_r \le 1$, which together with $m \ge 2$ gives $w_r' \in (0, 1)$ for $r = 1,\ldots,m$.
Applying the induction hypothesis to $U'$ with weights $w_1',\ldots,w_m'$ yields $I = \{i, j\} \subseteq \{1,\ldots,m\}$ ($i < j$) such that $U'_I = \diag(\sqrt{w_i'},\sqrt{w_j'}) A_I$ is invertible and $\|A_I^{-1}\|_2 = \bigl\|(U'_I)^{-1}\diag(\sqrt{w_i'},\sqrt{w_j'})\bigr\|_2 \le 1$.
\end{proof}

\section{Streaming row selection}
\label{sec:algorithm}

\Cref{sec:prototype} presents a high-level prototype algorithm that applies the elimination step of \Cref{sec:proof} to four rows at a time while retaining at most three rows between successive elimination steps.
\Cref{sec:closed-form} derives closed-form formulas for the local subproblems and develops a practical implementation.

\subsection{Prototype streaming algorithm}
\label{sec:prototype}

Let $U \in \R^{n \times 2}$ and $w_1,\ldots,w_n \in (0, 1)$ satisfy the assumptions of \Cref{thm:weighted}, and write $u_r^\top$ for the $r$-th row of $U$ and $a_r^\top = w_r^{-1/2} u_r^\top$ for the $r$-th row of $A = \diag(w_1^{-1/2},\ldots,w_n^{-1/2}) U$, so that \eqref{eq:weighted-bound} holds if and only if $A_I$ is invertible with $\|A_I^{-1}\|_2 \le 1$.
Here and throughout, we set $\|B^{-1}\|_2 = \infty$ whenever a square matrix $B$ is singular.
In the inductive step of \Cref{sec:proof}, we eliminated a row using a nonzero vector $\gamma \in \R^n$ satisfying $\sum_{r=1}^n \gamma_r a_r a_r^\top = 0$ and $\sum_{r=1}^n \gamma_r \le 0$.
Crucially, any four matrices $a_r a_r^\top$ already lie in the three-dimensional space of $2 \times 2$ real symmetric matrices and are linearly dependent, so we can eliminate a row using only four rows at a time while leaving all remaining rows untouched.

Zero rows are skipped as $U$ is read row by row since the final submatrix $U_I$ must be nonsingular.
If $U$ has at most three nonzero rows, no elimination is needed, and we directly evaluate $\|A_I^{-1}\|_2$ on the at most three $2 \times 2$ submatrices formed by the nonzero rows to find a pair satisfying $\|A_I^{-1}\|_2 \le 1$, which by \Cref{thm:weighted} is guaranteed to exist.
We may now therefore assume without loss of generality that $n \ge 4$ and all rows of $U$ are nonzero.
At row $4$, the four matrices $a_1 a_1^\top,\ldots,a_4 a_4^\top$ are linearly dependent, so there exists a nonzero vector $\gamma = (\gamma_1, \gamma_2, \gamma_3, \gamma_4)^\top \in \R^4$ satisfying $\sum_{k=1}^4 \gamma_k a_k a_k^\top = 0$ and $\sum_{k=1}^4 \gamma_k \le 0$.
Setting
\[
 \theta_* = \min_{\gamma_k < 0} \frac{w_k}{-\gamma_k} > 0, \qquad w_k^{(4)} = w_k + \theta_* \gamma_k \ge 0, \qquad k = 1,2,3,4,
\]
zeroes out at least one weight and leaves a subset $P_4 = \{p \in \{1,2,3,4\} : w_p^{(4)} > 0\}$ of size $|P_4| \le 3$ with positive weights $\{w_p^{(4)}\}_{p \in P_4}$ satisfying
\begin{equation}\label{eq:proto-r4}
 \sum_{p \in P_4} w_p^{(4)} a_p a_p^\top = \sum_{s=1}^4 u_s u_s^\top, \qquad \sum_{p \in P_4} w_p^{(4)} \le \sum_{s=1}^4 w_s,
\end{equation}
since $w_k a_k a_k^\top = u_k u_k^\top$.

For $n \ge 5$, we next add row $5$ with initial weight $w_5^{(4)} = w_5$ and rescaled vector $a_5 = w_5^{-1/2} u_5$, forming $P_4 \cup \{5\}$ of size at most four.
If $|P_4| \le 2$, then $|P_4 \cup \{5\}| \le 3$, so no elimination is needed and we set $P_5 = P_4 \cup \{5\}$ with $w_p^{(5)} = w_p^{(4)}$ for $p \in P_5$.
If $|P_4| = 3$, write $P_4 \cup \{5\} = \{p_1, p_2, p_3, p_4\}$ (with $p_4 = 5$) and compute a nonzero $\gamma = (\gamma_1, \gamma_2, \gamma_3, \gamma_4)^\top \in \R^4$ satisfying $\sum_{k=1}^4 \gamma_k a_{p_k} a_{p_k}^\top = 0$ and $\sum_{k=1}^4 \gamma_k \le 0$.
Updating the weights from $w_{p_k}^{(4)}$ via
\[
 \theta_* = \min_{\gamma_k < 0} \frac{w_{p_k}^{(4)}}{-\gamma_k} > 0, \qquad w_{p_k}^{(5)} = w_{p_k}^{(4)} + \theta_* \gamma_k \ge 0, \qquad k = 1,2,3,4,
\]
again zeroes out at least one weight and leaves $P_5 = \{p \in P_4 \cup \{5\} : w_p^{(5)} > 0\}$ of size $|P_5| \le 3$.
In either case, adding $w_5 a_5 a_5^\top = u_5 u_5^\top$ to \eqref{eq:proto-r4} gives
\[
 \sum_{p \in P_5} w_p^{(5)} a_p a_p^\top = \sum_{p \in P_4} w_p^{(4)} a_p a_p^\top + u_5 u_5^\top = \sum_{s=1}^5 u_s u_s^\top, \qquad
 \sum_{p \in P_5} w_p^{(5)} \le \sum_{p \in P_4} w_p^{(4)} + w_5 \le \sum_{s=1}^5 w_s.
\]

In general, repeating this step for $r = 4,\ldots,n$ yields a subset $P_r \subseteq \{1,\ldots,r\}$ of size $|P_r| \le 3$ with positive weights $\{w_p^{(r)}\}_{p \in P_r}$ satisfying
\[
 \sum_{p \in P_r} w_p^{(r)} a_p a_p^\top = \sum_{s=1}^r u_s u_s^\top, \qquad \sum_{p \in P_r} w_p^{(r)} \le \sum_{s=1}^r w_s.
\]
In particular, at $r = n$, the right-hand sides equal $U^\top U = I_2$ and $\sum_{s=1}^n w_s \le 1$.
Hence the $|P_n| \times 2$ matrix with rows $\sqrt{w_p^{(n)}}\,a_p^\top$ ($p \in P_n$) has orthonormal columns, which implies $2 \le |P_n| \le 3$ and $w_p^{(n)} \in (0, 1)$ for all $p \in P_n$.
Applying \Cref{thm:weighted} to this matrix guarantees a two-element subset $I \subseteq P_n$ with $\|A_I^{-1}\|_2 \le 1$, which is found by evaluating $\|A_I^{-1}\|_2$ on the at most three pairs in $P_n$.
\Cref{alg:proto} summarizes this prototype procedure.
Here, we update $P$ and $w$ in place.

\begin{algorithm}[t]
 \caption{Streaming $2 \times 2$ row selection (prototype)}
 \label{alg:proto}
 \begin{algorithmic}[1]
  \Require Matrix $U \in \R^{n \times 2}$ ($n \ge 2$) with $U^\top U = I_2$, weights $w \in (0, 1)^n$ with $\sum_{r=1}^n w_r \le 1$ (default $w_r = 1/n$)
  \Ensure Row index pair $I = \{i, j\}$ ($1 \le i < j \le n$) satisfying \eqref{eq:weighted-bound}
  \State $P \gets \emptyset$
  \For{$r = 1,\ldots,n$}
   \IIf{$u_r = 0$} \textbf{continue} \EndIIf
   \State $a_r \gets u_r / \sqrt{w_r}$, \quad $P \gets P \cup \{r\}$
   \If{$|P| = 4$} \Comment{write $P = \{p_1, p_2, p_3, p_4\}$}
    \State Find a nonzero $\gamma \in \R^4$ satisfying $\sum_{k=1}^4 \gamma_k a_{p_k} a_{p_k}^\top = 0$ and $\sum_{k=1}^4 \gamma_k \le 0$ \label{line:proto-null}
    \State $\theta_* \gets \min_{k : \gamma_k < 0} w_{p_k} / (-\gamma_k)$
    \State $w_{p_k} \gets w_{p_k} + \theta_* \gamma_k$ for $k = 1,2,3,4$, and remove all $p_k$ with $w_{p_k} = 0$ from $P$
   \EndIf
  \EndFor
  \State Return $I = \{i, j\} \subseteq P$ minimizing $\|A_I^{-1}\|_2$ \label{line:proto-final}
 \end{algorithmic}
\end{algorithm}

\subsection{Closed-form formulas and practical implementation}
\label{sec:closed-form}

Implementing \Cref{alg:proto} requires solving two local subproblems: finding a null vector $\gamma$ in \cref{line:proto-null} and evaluating $\|A_I^{-1}\|_2$ in \cref{line:proto-final}.
Rather than calling general-purpose SVD routines, which add overhead and can lose accuracy when the norms $\|a_{p_k}\|_2$ vary across many orders of magnitude, we derive numerically robust closed-form formulas for both subproblems.
We first consider the null-vector computation.

\begin{proposition}[Explicit formula for a null vector]\label{prop:nullvec}
Let $a_1, a_2, a_3, a_4 \in \R^2$ be nonzero vectors, and define
\[
 d_{pq} = a_{p,1}a_{q,2} - a_{p,2}a_{q,1}, \qquad 1 \le p, q \le 4.
\]
Then a nonzero vector $\gamma = (\gamma_1, \gamma_2, \gamma_3, \gamma_4)^\top \in \R^4$ satisfying $\sum_{q=1}^4 \gamma_q a_q a_q^\top = 0$ is given as follows:
\begin{enumerate}
 \item If two of the vectors are parallel, i.e., $d_{pq} = 0$ for some $1 \le p < q \le 4$, we may take
 \begin{equation}\label{eq:nullvec-par}
  \gamma = \|a_p\|_2^{-2} e_p - \|a_q\|_2^{-2} e_q.
 \end{equation}
 \item If no two of the vectors are parallel, i.e., $d_{pq} \ne 0$ for all $1 \le p < q \le 4$, we may take
 \begin{equation}\label{eq:nullvec}
  \gamma = \begin{bmatrix} d_{23}d_{24}d_{34} & -d_{13}d_{14}d_{34} & d_{12}d_{14}d_{24} & -d_{12}d_{13}d_{23} \end{bmatrix}^\top.
 \end{equation}
\end{enumerate}
\end{proposition}
\begin{proof}
If two vectors $a_p$ and $a_q$ ($1 \le p < q \le 4$) are parallel, then $a_p / \|a_p\|_2 = \pm a_q / \|a_q\|_2$ and $\|a_p\|_2^{-2} a_p a_p^\top - \|a_q\|_2^{-2} a_q a_q^\top = 0$, which gives \eqref{eq:nullvec-par}.

Suppose now that $d_{pq} \ne 0$ for all $1 \le p < q \le 4$.
Writing out $\sum_{q=1}^4 \gamma_q a_q a_q^\top = 0$ gives the $3 \times 4$ homogeneous linear system $M \gamma = 0$ with
\[
 M = \begin{bmatrix} m_1 & m_2 & m_3 & m_4 \end{bmatrix}
 = \begin{bmatrix}
  a_{1,1}^2 & a_{2,1}^2 & a_{3,1}^2 & a_{4,1}^2 \\
  a_{1,1}a_{1,2} & a_{2,1}a_{2,2} & a_{3,1}a_{3,2} & a_{4,1}a_{4,2} \\
  a_{1,2}^2 & a_{2,2}^2 & a_{3,2}^2 & a_{4,2}^2
 \end{bmatrix} \in \R^{3 \times 4}.
\]
Then
\[
 \gamma = \begin{bmatrix} \det[m_2, m_3, m_4] & -\det[m_1, m_3, m_4] & \det[m_1, m_2, m_4] & -\det[m_1, m_2, m_3] \end{bmatrix}^\top
\]
is a solution to $M \gamma = 0$.
Indeed, the $i$-th entry of $M \gamma$ is the determinant of the $4 \times 4$ matrix obtained by appending the $i$-th row of $M$ to the top of $M$, i.e.,
\[
 (M\gamma)_i = e_i^\top M \gamma = \det\begin{bmatrix} e_i^\top M \\ M \end{bmatrix} = 0, \qquad i = 1, 2, 3,
\]
and hence $M \gamma = 0$.
Now we compute the determinants $\det[m_p, m_q, m_r]$ for $1 \le p < q < r \le 4$ in terms of the $d_{pq}$.
Note that, for $a_{p,2}a_{q,2}a_{r,2} \ne 0$, factoring $a_{l,2}^2$ from column $m_l$ ($l \in \{p, q, r\}$) and setting $t_l = a_{l,1} / a_{l,2}$ gives a Vandermonde determinant:
\[
 \det[m_p, m_q, m_r]
 = a_{p,2}^2 a_{q,2}^2 a_{r,2}^2 \det\begin{bmatrix} t_p^2 & t_q^2 & t_r^2 \\ t_p & t_q & t_r \\ 1 & 1 & 1 \end{bmatrix}
 = a_{p,2}^2 a_{q,2}^2 a_{r,2}^2 (t_p - t_q)(t_p - t_r)(t_q - t_r)
 = d_{pq} d_{pr} d_{qr}.
\]
For $a_{p,2}a_{q,2}a_{r,2} = 0$, one can check that the same identity holds.
Substituting these determinants into the formula for $\gamma$ gives \eqref{eq:nullvec}.
\end{proof}

\begin{algorithm}[H]
 \caption{Streaming $2 \times 2$ row selection}
 \label{alg:opt2row}
 \begin{algorithmic}[1]
  \Require Matrix $U \in \R^{n \times 2}$ ($n \ge 2$) with $U^\top U = I_2$, weights $w \in (0, 1)^n$ with $\sum_{r=1}^n w_r \le 1$ (default $w_r = 1/n$)
  \Ensure Row index pair $I = \{i, j\}$ ($1 \le i < j \le n$) satisfying \eqref{eq:weighted-bound}
  \State $P \gets \emptyset$
  \For{$r = 1,\ldots,n$} \Comment{single streaming pass}
   \IIf{$u_r = 0$} \textbf{continue} \EndIIf \Comment{skip zero row}
   \State $a_r \gets u_r / \sqrt{w_r}$, \quad $\hat a_r \gets a_r / \|a_r\|_2$, \quad $P \gets P \cup \{r\}$
   \If{$|P| = 4$} \Comment{eliminate at least one row from $P = \{p_1, p_2, p_3, p_4\}$}
    \State $d_{kl} \gets \hat a_{p_k,1}\hat a_{p_l,2} - \hat a_{p_k,2}\hat a_{p_l,1}$ for $1 \le k < l \le 4$
    \If{$|d_{kl}| \le 4u$ for some $k < l$}
     \State $\gamma \gets \|a_{p_k}\|_2^{-2} e_k - \|a_{p_l}\|_2^{-2} e_l$ \Comment{\eqref{eq:nullvec-par}}
    \Else
     \State $\gamma \gets [d_{23}d_{24}d_{34},\, -d_{13}d_{14}d_{34},\, d_{12}d_{14}d_{24},\, -d_{12}d_{13}d_{23}]^\top$ \Comment{\eqref{eq:nullvec}}
     \State $\gamma_k \gets \gamma_k \|a_{p_k}\|_2^{-2}$ for $k = 1,2,3,4$
    \EndIf
    \State $\gamma \gets \pm \gamma / \|\gamma\|_\infty$ such that $\gamma_1 + \gamma_2 + \gamma_3 + \gamma_4 \le 0$
    \State Choose $k_* \in \arg\min_{k : \gamma_k < 0} w_{p_k} / (-\gamma_k)$ and set $\theta_* \gets w_{p_{k_*}} / (-\gamma_{k_*})$
    \State $w_{p_k} \gets w_{p_k} + \theta_* \gamma_k$ for $k \ne k_*$, and remove $p_{k_*}$ and all $p_k$ with $w_{p_k} \le 0$ from $P$
   \EndIf
  \EndFor
  \State Return $I = \{i, j\} \subseteq P$ minimizing $\|A_I^{-1}\|_2$ evaluated via \eqref{eq:inv2norm} \label{line:final}
 \end{algorithmic}
\end{algorithm}

Evaluating \eqref{eq:nullvec} directly from $a_1,\ldots,a_4$ can be numerically unstable in floating-point arithmetic: the determinants $d_{pq}$ scale with $\|a_p\|_2\|a_q\|_2$, so the products of three determinants in \eqref{eq:nullvec} have degree six in the norms $\|a_q\|_2$ and can overflow or underflow at moderate scales.
Furthermore, $|d_{pq}|$ cannot distinguish nearly parallel vectors from vectors of small length.
Both issues are mitigated by normalizing each vector to $\hat a_q = a_q / \|a_q\|_2 \in \R^2$.
Then we compute the determinants $\hat d_{pq} = \hat a_{p,1}\hat a_{q,2} - \hat a_{p,2}\hat a_{q,1} \in [-1, 1]$ from the normalized vectors.
When two vectors are parallel, i.e., numerically we check $|\hat d_{pq}| \le 4u$ where $u$ is the unit roundoff, we evaluate $\gamma$ directly by \eqref{eq:nullvec-par}.
In the nonparallel case, letting $\hat\gamma \in \R^4$ denote \eqref{eq:nullvec} evaluated with $\hat d_{pq}$ in place of $d_{pq}$, we have
\[
 \gamma_q = \Bigl(\prod_{l=1}^4 \|a_l\|_2^2\Bigr)\,\|a_q\|_2^{-2}\,\hat\gamma_q, \qquad q = 1,2,3,4.
\]
Since any nonzero scalar multiple of a null vector remains a null vector, we further drop the common factor $\prod_{l=1}^4 \|a_l\|_2^2$.
Finally, in both cases, we normalize the resulting null vector in the $\infty$-norm.

For the final $\|A_I^{-1}\|_2$ check in \cref{line:proto-final} of \Cref{alg:proto}, the next proposition gives a closed-form formula.

\begin{proposition}[Explicit formula for $\|B^{-1}\|_2$ when $B\in\R^{2 \times 2}$]\label{prop:inv2norm}
Let $B = \begin{bsmallmatrix} b_{11} & b_{12} \\ b_{21} & b_{22} \end{bsmallmatrix} \in \R^{2 \times 2}$ have singular values $\sigma_1(B) \ge \sigma_2(B) \ge 0$.
Then
\[
 \sigma_{1,2}(B) = \frac{\sqrt{\|B\|_\F^2 + 2|\det B|} \pm \sqrt{\|B\|_\F^2 - 2|\det B|}}{2},
\]
and when $B$ is nonsingular,
\begin{equation}\label{eq:inv2norm}
 \|B^{-1}\|_2 = \frac{1}{\sigma_2(B)} = \frac{\sigma_1(B)}{|\det B|}
 = \frac{\sqrt{(b_{11}+b_{22})^2 + (b_{12}-b_{21})^2}+\sqrt{(b_{11}-b_{22})^2 + (b_{12}+b_{21})^2}}{2\,|b_{11}b_{22} - b_{12}b_{21}|}.
\end{equation}
\end{proposition}
\begin{proof}
Write $\sigma_i = \sigma_i(B)$ for $i = 1, 2$.
Since $\sigma_1^2 + \sigma_2^2 = \|B\|_\F^2$ and $\sigma_1\sigma_2 = |\det B|$, we have $(\sigma_1 \pm \sigma_2)^2 = \|B\|_\F^2 \pm 2|\det B|$.
Because $\sigma_1 \ge \sigma_2 \ge 0$, taking square roots and adding or subtracting the two equations gives the formula for $\sigma_{1,2}(B)$.
Explicit computation gives $\|B\|_\F^2 \pm 2\det B = (b_{11} \pm b_{22})^2 + (b_{12} \mp b_{21})^2$, which leads to the desired formula for $\|B^{-1}\|_2$.
\end{proof}

The rationalized form \eqref{eq:inv2norm} avoids subtracting the two square roots in $\sigma_2(B)$.
\Cref{alg:opt2row} summarizes the resulting practical implementation.
Each iteration of \Cref{alg:opt2row} performs $\mathcal O(1)$ arithmetic operations on at most four vectors in $\R^2$, and the final check in \cref{line:final} evaluates \eqref{eq:inv2norm} on at most $\binom{3}{2} = 3$ pairs.
In total, \Cref{alg:opt2row} requires $\mathcal O(n)$ operations, a single sequential read of $U$ and $w$, and $\mathcal O(1)$ auxiliary memory.

\section{Counterexamples for standard selection strategies}
\label{sec:counterexamples}

We next show that standard volume-based and greedy pivoting strategies fail to attain the GTZ bound $\sqrt{n}$ in \eqref{eq:gtz}, even for $k = 2$.
Specifically, for $U \in \R^{n \times k}$ with $U^\top U = I_k$, we consider five standard strategies (writing $I_0 = \emptyset$ and $I_t = I_{t-1} \cup \{j_t\}$ for $t = 1,\ldots,k$ in the three sequential schemes):
\begin{itemize}
 \item \emph{Global maximum-volume selection} (exact maxvol~\cite{Knuth1985,GTZ1997theory,GTZ1997maxvol,Goreinov2010GoodSubmatrix}), which chooses $|I| = k$ maximizing $|\det U_I|$ over all $\binom{n}{k}$ subsets.
 \item \emph{Local maximum-volume swapping} (maxvol~\cite{Goreinov2010GoodSubmatrix}), which starts from a nonsingular $U_I$ (typically obtained by LUPP or QRCP) and repeatedly replaces a row in $I$ whenever doing so increases $|\det U_I|$ by more than a factor $1 + \varepsilon$ ($\varepsilon \ge 0$), terminating once $\|U U_I^{-1}\|_{\max} \le 1 + \varepsilon$~\cite{Knuth1985,GuEisenstat1996,Goreinov2010GoodSubmatrix}.
 \item \emph{Column-pivoted QR on $U^\top$} (QRCP~\cite{BusingerGolub1965}, used in QDEIM~\cite{DrmacGugercin2016}), which at step $t$ selects (cf.~\cite[Remark~2.3]{DrmacGugercin2016})
 \[
  j_t \in \arg\max_{r \notin I_{t-1}} \vol\bigl(U_{I_{t-1}\cup\{r\}}\bigr), \qquad \vol(M) = \sqrt{\det(M M^\top)},
 \]
 so that at $k = 2$, $j_1$ maximizes $\|u_r\|_2$ and $j_2 \ne j_1$ maximizes $|\det U_{\{j_1,r\}}|$.
 \item \emph{Frobenius-norm greedy pivoting} (BSQR~\cite{Bischof1990,Stewart1990,Osinsky2025,Damle2026}), which at step $t$ selects (cf.~\cite[(2.3)--(2.4)]{Damle2026})
 \[
  j_t \in \arg\min_{\substack{r \notin I_{t-1} \\ \rank(U_{I_{t-1}\cup\{r\}}) = t}} \bigl\|U_{I_{t-1}\cup\{r\}}^\dagger\bigr\|_\F^2,
 \]
 so that at $k = 2$, $j_1$ maximizes $\|u_r\|_2$ and $j_2$ minimizes $(\|u_{j_1}\|_2^2 + \|u_r\|_2^2)/(\det U_{\{j_1,r\}})^2$.
 \item \emph{LU with partial pivoting} (LUPP, equivalent to DEIM~\cite{ChaturantabutSorensen2010}), which at step $t$ selects (cf.~\cite[Section~2.1.2]{DrmacGugercin2016})
 \[
  j_t \in \arg\max_{r \notin I_{t-1}} \bigl|\det U_{I_{t-1}\cup\{r\},\,1{:}t}\bigr|,
 \]
 so that at $k = 2$, $j_1$ maximizes $|u_{r,1}|$ and $j_2 \ne j_1$ maximizes $|\det U_{\{j_1,r\}}|$.
\end{itemize}
Exact maxvol, maxvol with $\varepsilon = 0$, and BSQR guarantee $\|U_I^{-1}\|_2 \le \sqrt{k(n-k)+1}$~\cite{GuEisenstat1996,GTZ1997theory,Goreinov2010GoodSubmatrix,Osinsky2025,Damle2026} (which is sharp for volume maximization~\cite{GTZ1997theory}), whereas the worst-case bounds for QRCP and LUPP grow exponentially with $k$ on Kahan-type matrices~\cite{GuEisenstat1996,ChaturantabutSorensen2010,DrmacGugercin2016}.
At $k = 2$, the bound $\sqrt{2n - 3}$ coincides with $\sqrt{n}$ when $n = 3$, and strictly exceeds $\sqrt{n}$ for every $n \ge 4$.

To construct counterexamples, we consider a family of three-block orthonormal bases generalizing the equality cases of \eqref{eq:gtz} for $k = 2$ characterized in~\cite{Nesterenko2026equality}.
Let $n \ge 3$ and let $p_1, p_2, p_3$ be positive integers satisfying $p_1 + p_2 + p_3 = n$.
Given a unit vector $z = (z_1, z_2, z_3)^\top \in \R^3$ and an orthonormal basis $V = [v_1, v_2, v_3]^\top \in \R^{3 \times 2}$ of $z^\perp$, i.e., $V^\top V = I_2$ and $V^\top z = 0$, define
\begin{equation}\label{eq:U-block}
 U_{p_1,p_2,p_3}(z) = \begin{bmatrix} p_1^{-1/2}\mathbf{1}_{p_1} v_1^\top \\[2pt] p_2^{-1/2}\mathbf{1}_{p_2} v_2^\top \\[2pt] p_3^{-1/2}\mathbf{1}_{p_3} v_3^\top \end{bmatrix} \in \R^{n \times 2},
\end{equation}
where $\mathbf{1}_m \in \R^m$ denotes the all-ones vector.
Direct multiplication gives $U_{p_1,p_2,p_3}(z)^\top U_{p_1,p_2,p_3}(z) = \sum_{i=1}^3 p_i(p_i^{-1} v_i v_i^\top) = V^\top V = I_2$.

\begin{theorem}[Equality cases for $k = 2$, cf.~{\cite[Proposition~1]{Nesterenko2026equality}}]\label{thm:equality}
Let $n \ge 3$ and let $U \in \R^{n \times 2}$ satisfy $U^\top U = I_2$.
Then
\begin{equation}\label{eq:equality-iff}
 \min_{|I|=2}\|U_I^{-1}\|_2 = \sqrt{n} \iff U = D P U_{p_1,p_2,p_3}(z_*)
\end{equation}
for some positive integers $p_1 + p_2 + p_3 = n$, permutation matrix $P \in \R^{n \times n}$, diagonal matrix $D$ with $\pm 1$ on the diagonal, and orthonormal basis $V \in \R^{3 \times 2}$ of $z_*^\perp$, where
\begin{equation}\label{eq:z-star}
 z_* = \left(\sqrt{\frac{n - p_1}{2n}},\, \sqrt{\frac{n - p_2}{2n}},\, \sqrt{\frac{n - p_3}{2n}}\right)^\top.
\end{equation}
Moreover, every $2 \times 2$ submatrix $B$ of $U_{p_1,p_2,p_3}(z_*)$ with rows from distinct blocks satisfies $\|B^{-1}\|_2 = \sqrt{n}$.
\end{theorem}
\begin{proof}
For any unit vector $z \in \R^3$, any $2 \times 2$ submatrix of $U_{p_1,p_2,p_3}(z)$ taking two rows from the same block is singular.
Let $B \in \R^{2 \times 2}$ take one row from block $i$ and one row from block $j$ ($1 \le i < j \le 3$).
To compute $\|B^{-1}\|_2$ via \Cref{prop:inv2norm}, note from \eqref{eq:U-block} that $B = \begin{bsmallmatrix} p_i^{-1/2} v_i^\top \\ p_j^{-1/2} v_j^\top \end{bsmallmatrix}$, so
\begin{equation}\label{eq:block-invariants}
 \|B\|_\F^2 = \frac{\|v_i\|_2^2}{p_i} + \frac{\|v_j\|_2^2}{p_j}, \qquad (\det B)^2 = \frac{(\det[v_i, v_j])^2}{p_i p_j}.
\end{equation}
Since $[V, z] \in \R^{3 \times 3}$ is orthogonal, $V V^\top = I_3 - z z^\top$ gives $[v_i, v_j]^\top [v_i, v_j] = \begin{bsmallmatrix} 1 - z_i^2 & -z_i z_j \\ -z_i z_j & 1 - z_j^2 \end{bsmallmatrix}$, whose diagonal entries and determinant yield $\|v_i\|_2^2 = 1 - z_i^2$ and $(\det[v_i, v_j])^2 = 1 - z_i^2 - z_j^2$.
For $z_*$ in \eqref{eq:z-star}, we have $\|v_i\|_2^2 = \frac{n + p_i}{2n}$ and $(\det[v_i, v_j])^2 = \frac{p_i + p_j}{2n}$.
Substituting them into \eqref{eq:block-invariants} yields
\[
 \|B\|_\F^2 = \frac{1}{n} + \frac{1}{2}\Bigl(\frac{1}{p_i} + \frac{1}{p_j}\Bigr), \qquad |\det B| = \sqrt{\frac{1}{2n}\Bigl(\frac{1}{p_i} + \frac{1}{p_j}\Bigr)}.
\]
Since $\frac{1}{2}\bigl(\frac{1}{p_i} + \frac{1}{p_j}\bigr) \ge \frac{1}{p_i + p_j} > \frac{1}{n}$, we have $\sqrt{\|B\|_\F^2 \pm 2|\det B|} = \sqrt{\frac{1}{2}\bigl(\frac{1}{p_i} + \frac{1}{p_j}\bigr)} \pm \frac{1}{\sqrt{n}}$, so \Cref{prop:inv2norm} gives $\|B^{-1}\|_2 = \sqrt{n}$.
Because permuting rows via $P$ and scaling rows by $\pm 1$ via $D$ preserve the singular values of all $2 \times 2$ row submatrices, every $U = D P U_{p_1,p_2,p_3}(z_*)$ satisfies $\min_{|I|=2}\|U_I^{-1}\|_2 = \sqrt{n}$.

Conversely, suppose $\min_{|I|=2}\|U_I^{-1}\|_2 = \sqrt{n}$, write $u_r^\top = (u_{r,1}, u_{r,2})$ for the $r$-th row of $U$, and set $q_r = (u_{r,1}^2 - u_{r,2}^2,\, 2u_{r,1}u_{r,2})^\top \in \R^2$ for $r = 1,\ldots,n$.
By~\cite[Proposition~1]{Nesterenko2026equality}, $\{1,\ldots,n\}$ partitions into three nonempty subsets $C_1, C_2, C_3$ of sizes $p_1 + p_2 + p_3 = n$ such that $q_r$ is constant on each $C_i$ and satisfies $\|q_r\|_2 = \frac{1}{2n} + \frac{1}{2p_i}$ ($i = 1, 2, 3$).
Note that $\|q_r\|_2 = \|u_r\|_2^2$ and that $q_r = q_s$ if and only if $u_r = \pm u_s$.
On each $C_i$, we may therefore write $u_r = \pm p_i^{-1/2} v_i$ for some vector $v_i \in \R^2$ satisfying $\|v_i\|_2^2 = p_i\|u_r\|_2^2 = \frac{n + p_i}{2n}$.
Setting $V = [v_1, v_2, v_3]^\top \in \R^{3 \times 2}$ gives $V^\top V = \sum_{i=1}^3 p_i(p_i^{-1} v_i v_i^\top) = \sum_{r=1}^n u_r u_r^\top = U^\top U = I_2$.
Let $z = (z_1, z_2, z_3)^\top \in \R^3$ be a unit vector satisfying $V^\top z = 0$.
By the orthogonality of $[V, z]$, we have $\|v_i\|_2^2 = 1 - z_i^2$, so that $z_i^2 = 1 - \|v_i\|_2^2 = \frac{n - p_i}{2n} = z_{*,i}^2$ for $i = 1, 2, 3$.
By changing the signs of $(v_i, z_i)$ if necessary, we may assume $z = z_*$, so $V$ is an orthonormal basis of $z_*^\perp$.
Permuting the blocks of $U_{p_1,p_2,p_3}(z_*)$ onto $C_1, C_2, C_3$ via $P$ and choosing the row signs in $D$ then yields $U = D P U_{p_1,p_2,p_3}(z_*)$.
\end{proof}

The equality characterization in \Cref{thm:equality} is closely related to the weighted three-row case in the proof of \Cref{thm:weighted}.
Indeed, starting with uniform weights $1/n$, merging identical rows within each block of $U_{p_1,p_2,p_3}(z)$ reduces the unweighted problem to a three-row weighted problem with orthonormal basis $V \in \R^{3 \times 2}$ and weights $w_i = p_i / n$.
For the unit vector $z \in \range(V)^\perp$, the key inequality \eqref{eq:n3-key} in \Cref{sec:proof} is strict unless all three ratios $z_i^2 / (1 - w_i)$ coincide.
Since $\sum_{i=1}^3 w_i = 1$ and $\|z\|_2 = 1$, this equality condition is equivalent to $z_i^2 = (1 - w_i) / 2 = (n - p_i) / (2n)$, precisely recovering $z_*$ in \eqref{eq:z-star} up to signs.
Thus, the unweighted $n$-row equality cases are obtained simply by repeating and rescaling the rows of the weighted three-row equality cases.

Notice also from the proof of \Cref{thm:equality} that on $U_{p_1,p_2,p_3}(z_*)$, every $2 \times 2$ submatrix $B$ with rows from distinct blocks has the same minimum singular value $\sigma_2(B) = n^{-1/2}$, whereas $|\det B| = \sqrt{(p_i + p_j)/(2n p_i p_j)}$ varies with the block sizes.
Choosing unbalanced block sizes $(p_1, p_2, p_3) = (1, m, m)$ ($n = 2m + 1$) gives pairs $\{1, s\}$ ($s \ge 2$) a much larger volume at $z_*$ than pairs between the two large blocks, suggesting a counterexample for volume maximization upon varying $z$ away from $z_*$.
It turns out that a suitable choice of $z$ and $V$ makes all five standard strategies select submatrices whose inverse spectral norms exceed $\sqrt{n}$ by the same asymptotic constant factor.

\begin{proposition}[Simultaneous suboptimality of standard strategies]\label{prop:suboptimal}
For every odd $n = 2m + 1 \ge 9$, let $U = U_{1,m,m}(z) \in \R^{n \times 2}$ in \eqref{eq:U-block} be defined with $z = (\sqrt{1 - 4/n}, \sqrt{2/n}, \sqrt{2/n})^\top$ and an orthonormal basis $V \in \R^{3 \times 2}$ of $z^\perp$ satisfying $v_{1,2} = 0$.
For every even $n \ge 12$, let $U \in \R^{n \times 2}$ be obtained by appending a zero row to the $(n - 1) \times 2$ matrix just defined.
Then exact maxvol, maxvol (with $\varepsilon = 0$ from any nonsingular initial pair, or with any $\varepsilon \ge 0$ when initialized by LUPP or QRCP), QRCP, BSQR, and LUPP all select a submatrix $U_I$, whereas \Cref{alg:opt2row} (with $w_r = 1/n$) selects a submatrix $U_J$ attaining $\min_{|K|=2}\|U_K^{-1}\|_2$, such that
\begin{equation}\label{eq:suboptimal-gap}
 \|U_I^{-1}\|_2 > \sqrt{n} > \|U_J^{-1}\|_2, \qquad \lim_{n \to \infty} \frac{\|U_I^{-1}\|_2}{\sqrt{n}} = \frac{\sqrt{3 + \sqrt{5}}}{2}, \qquad \lim_{n \to \infty} \frac{\|U_I^{-1}\|_2}{\|U_J^{-1}\|_2} = \frac{1 + \sqrt{5}}{2}.
\end{equation}
\end{proposition}
\begin{proof}
Suppose first that $n = 2m + 1 \ge 9$ is odd, so $p_1 = 1$ and $p_2 = p_3 = m = \frac{n - 1}{2}$.
A valid basis of $z^\perp$ with $v_{1,2} = 0$ is $V = \begin{bsmallmatrix} 2/\sqrt{n} & 0 \\ -\sqrt{(n-4)/(2n)} & 1/\sqrt{2} \\ -\sqrt{(n-4)/(2n)} & -1/\sqrt{2} \end{bsmallmatrix}$.
We have $\|v_1\|_2^2 = \frac{4}{n}$ and $\|v_2\|_2^2 = \|v_3\|_2^2 = \frac{n - 2}{n}$, so row $1$ of $U$ has $\|u_1\|_2^2 = \frac{4}{n}$ whereas every row $s \ge 2$ has $\|u_s\|_2^2 = \frac{n - 2}{m n} = \frac{2(n - 2)}{n(n - 1)} < \frac{4}{n} = \|u_1\|_2^2$.
For every pair $I = \{1, s\}$ ($s \ge 2$) and every pair $K \subseteq \{2,\ldots,n\}$ between blocks $2$ and $3$, substituting $(\det[v_1, v_j])^2 = 1 - z_1^2 - z_j^2 = \frac{2}{n}$ ($j = 2, 3$) and $(\det[v_2, v_3])^2 = 1 - z_2^2 - z_3^2 = \frac{n - 4}{n}$ into \eqref{eq:block-invariants} gives
\[
 \|U_I\|_\F^2 = \frac{6n - 8}{n(n - 1)}, \quad |\det U_I| = \frac{2}{\sqrt{n(n - 1)}}, \quad \|U_K\|_\F^2 = \frac{4(n - 2)}{n(n - 1)}, \quad |\det U_K| = \frac{2\sqrt{n - 4}}{(n - 1)\sqrt{n}}.
\]
Substituting these quantities into \Cref{prop:inv2norm} yields, for all $n \ge 9$,
\[
 \|U_I^{-1}\|_2 = \frac{1}{2}\sqrt{3n - 4 + \sqrt{5n^2 - 20n + 16}} > \sqrt{n}, \qquad \|U_K^{-1}\|_2 = \sqrt{\frac{n(n - 1)}{2(n - 4)}} < \sqrt{n}.
\]
Because $v_{1,2} = 0$ gives $|u_{1,1}| = \|u_1\|_2 > \|u_s\|_2 \ge |u_{s,1}|$ for all $s \ge 2$, QRCP, BSQR, and LUPP uniquely select $j_1 = 1$ at step $1$.
Since $\|u_s\|_2^2$ and $|\det U_I| > 0$ are identical for all $I = \{1, s\}$ ($s \ge 2$), each of these three methods selects such a pair $I$ at step $2$.
Moreover, $n - 1 > n - 4$ gives $|\det U_I| > |\det U_K|$, so the global maximum-volume pairs are precisely $I = \{1, s\}$ ($s \ge 2$), and Cramer's rule gives $\|U U_I^{-1}\|_{\max} = 1$.
Hence maxvol terminates immediately at $I$ for any $\varepsilon \ge 0$ under LUPP or QRCP initialization, and when $\varepsilon = 0$ from any nonsingular pair $K \subseteq \{2,\ldots,n\}$ it swaps in row $1$ to reach such a pair $I$ in one step and terminates.
Thus all five standard strategies select a pair $I = \{1, s\}$ ($s \ge 2$).
In contrast, because \Cref{alg:opt2row} guarantees $\|U_J^{-1}\|_2 \le \sqrt{n}$ while every other nonsingular pair has the form $I = \{1, s\}$ and exceeds $\sqrt{n}$, $J$ must be a pair $K$ between blocks $2$ and $3$, attaining $\min_{|K|=2}\|U_K^{-1}\|_2$.
Dividing $\|U_I^{-1}\|_2$ and $\|U_J^{-1}\|_2$ by $\sqrt{n}$ and letting $n \to \infty$ proves the two limits in \eqref{eq:suboptimal-gap}.

For even $n \ge 12$, appending $u_n^\top = 0$ to the $(n - 1) \times 2$ matrix above preserves $U^\top U = I_2$, and since $\|u_n\|_2 = |u_{n,1}| = 0$ and every pair containing row $n$ is singular, all five standard strategies still select a pair $I = \{1, s\}$ ($2 \le s \le n - 1$).
Replacing $n$ by $n - 1$ in the odd-case formulas gives $\|U_I^{-1}\|_2 > \sqrt{n}$ for all $n \ge 12$, while every pair $K \subseteq \{2,\ldots,n-1\}$ between blocks $2$ and $3$ still satisfies $\|U_K^{-1}\|_2 < \sqrt{n - 1} < \sqrt{n}$.
Hence \Cref{alg:opt2row} (which skips $u_n = 0$ and guarantees $\|U_J^{-1}\|_2 \le \sqrt{n}$) again selects such a pair $J$ attaining $\min_{|K|=2}\|U_K^{-1}\|_2$, and $\sqrt{n - 1}/\sqrt{n} \to 1$ preserves both limits in \eqref{eq:suboptimal-gap}.
\end{proof}

\section{Numerical experiments}
\label{sec:experiments}

We compare \Cref{alg:opt2row} (with uniform weights $w_r = 1/n$) in MATLAB against the five standard strategies from \Cref{sec:counterexamples} (exact maxvol, maxvol, QRCP, BSQR, and LUPP), alongside $\min_{|I|=2}\|U_I^{-1}\|_2$ (labeled \emph{exact minimum}).
Here, exact maxvol and the exact minimum are computed by exhaustive search over all $\binom{n}{2}$ row pairs, QRCP and LUPP use MATLAB's built-in pivoted factorizations, and maxvol (initialized by LUPP with default swap tolerance $\varepsilon = 5 \cdot 10^{-2}$~\cite{Goreinov2010GoodSubmatrix}) and BSQR~\cite{Damle2026} are adapted from existing public implementations (see the \hyperref[sec:reproducibility]{Reproducibility Statement}).

\begin{figure}[t]
 \centering
 \includegraphics[width=\linewidth]{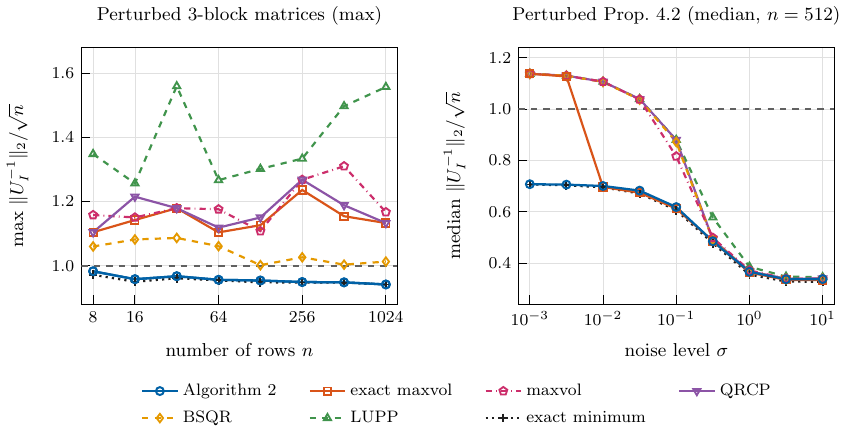}
 \caption{Left: maximum of $\|U_I^{-1}\|_2 / \sqrt{n}$ over $2{,}000$ perturbed three-block matrices $U$ in \eqref{eq:perturb} for $8 \le n \le 1024$ ($\sigma = 0.02$).
 Right: median of $\|U_I^{-1}\|_2 / \sqrt{n}$ over $2{,}000$ perturbations \eqref{eq:perturb} of the counterexample from \Cref{prop:suboptimal} at $n = 512$ across noise levels $10^{-3} \le \sigma \le 10^{1}$.
 The dashed horizontal line marks the sharp bound $1$.}
 \label{fig:experiments}
\end{figure}

We perform two experiments, summarized in \Cref{fig:experiments}.
First, for each $n \in \{8, 16, \ldots, 1024\}$, we generate $2{,}000$ matrices $U \in \R^{n \times 2}$ via the thin QR factorization
\begin{equation}\label{eq:perturb}
 U R = D P U_{p_1,p_2,p_3}(z) + \frac{\sigma}{\sqrt{n}}\,G
\end{equation}
at noise level $\sigma = 0.02$, where $D$ has random $\pm 1$ diagonal entries, $P$ is a random permutation matrix, and $G \in \R^{n \times 2}$ has i.i.d.\ standard normal entries.
The vector $z \in \R^3$ is drawn uniformly on the unit sphere, and we take a random orthonormal basis $V$ of $z^\perp$ in \eqref{eq:U-block}.
For the block sizes $p_1, p_2, p_3$, we draw positive integers $p_1 + p_2 + p_3 = n$ roughly proportional to $n^{u_i}$ with independent $u_i \sim \operatorname{Unif}(0, 1)$ ($i = 1, 2, 3$), which produces both balanced and strongly unbalanced block sizes as in \Cref{prop:suboptimal}.
The maximum of $\|U_I^{-1}\|_2 / \sqrt{n}$ over the $2{,}000$ trials is shown in \Cref{fig:experiments} (left).
Second, fixing $n = 512$, we generate $2{,}000$ matrices $U$ at each noise level $\sigma \in [10^{-3}, 10^{1}]$ in the same way via \eqref{eq:perturb} with $U_{p_1,p_2,p_3}(z)$ being the counterexample matrix of \Cref{prop:suboptimal}, and plot the median of $\|U_I^{-1}\|_2 / \sqrt{n}$ over the $2{,}000$ trials in \Cref{fig:experiments} (right).

In the first experiment (\Cref{fig:experiments}, left), all five standard strategies exceed $\sqrt{n}$ across all tested dimensions, whereas \Cref{alg:opt2row} remains strictly below $\sqrt{n}$ on every instance.
In the second experiment (\Cref{fig:experiments}, right), the standard strategies stay near the theoretical plateau $\frac{1}{2}\sqrt{3 + \sqrt{5}}$ above $\sqrt{n}$ for small perturbations ($\sigma \ll 1$), with exact maxvol dropping below $\sqrt{n}$ earlier than the remaining four strategies as $\sigma$ increases, before all methods descend well below $\sqrt{n}$ in the isotropic Gaussian regime ($\sigma \gtrsim 1$).
Across both experiments, \Cref{alg:opt2row} closely tracks the exact minimum.

\section*{Reproducibility Statement}
\label{sec:reproducibility}
\addcontentsline{toc}{section}{Reproducibility Statement}
All numerical experiments were run in MATLAB R2025a.
The source code and data reproducing \Cref{fig:experiments} are available at \url{https://github.com/bonans/optindsel}, which includes implementations of \Cref{alg:opt2row}, the three-block generator for \eqref{eq:U-block}, exhaustive search for exact maxvol and the exact minimum, and adaptations of maxvol from TT-Toolbox (\url{https://github.com/oseledets/TT-Toolbox}) and BSQR from \url{https://github.com/asdamle/BischofStewartQR}.

\section*{AI Usage Statement and Acknowledgments}
\addcontentsline{toc}{section}{AI Usage Statement and Acknowledgments}
Broader exploration of higher-rank cases of the GTZ conjecture made use of OpenAI GPT, Anthropic Claude, and Google Gemini models in multi-agent workflows.
The mathematical results presented in this paper were developed through human--AI interaction using Google Gemini 4 Carbon within Google Antigravity.
The weighted row-elimination idea emerged during an exploration of the $n \times 3$ case with Gemini.
The author identified its specialization to $k = 2$ and directed the development of a short proof.
The author then posed the question of turning this proof into a constructive row-selection algorithm and initiated the counterexample study against maximum-volume and greedy pivoting strategies.
Gemini also assisted with the implementation of the numerical experiments and the preparation of figures.
The manuscript was written by the author, with Gemini and OpenAI ChatGPT (GPT-6 and GPT-6.2) assisting in language editing and proofreading.
The author verified all results and assumes full responsibility for the manuscript.

The author gratefully acknowledges Google for providing access to Gemini models through the Gemini for Science program (\url{https://ai.google/gemini-for-science/}).

\addcontentsline{toc}{section}{References}
\bibliographystyle{abbrvurl}
\bibliography{references}
\end{document}